\documentclass[12pt,reqno]{amsart}

\usepackage[arrow,matrix,curve]{xy}

\usepackage[dvips]{graphicx} 

\usepackage{amssymb, latexsym, amsmath, amscd, array, epigraph, url,
setspace 
}

\usepackage
[breaklinks=true,
hypertexnames=false]   
{hyperref}

\usepackage{breakurl}   

\theoremstyle{definition} 

\numberwithin{equation}{section}

\author[J. B.]{Jacques Bair} \address{J. Bair, HEC-ULG, University
of Liege, 4000 Belgium} \email{j.bair@ulg.ac.be}

\author[P. B.]{Piotr B\l{}aszczyk} \address{P. B\l{}aszczyk, Institute
of Mathematics, Pedagogical University of Cracow, Poland}
\email{pb@up.krakow.pl}

\author[P. H.]{Peter Heinig} \address{P. Heinig}
\email{heinig@ma.tum.de}

\author[M. K.]{Mikhail G. Katz} \address{M. Katz, Department of
Mathematics, Bar Ilan University, Ramat Gan 52900 Israel}
\email{katzmik@macs.biu.ac.il}

\author[J. S.]{Jan Peter Sch\"afermeyer}
\address{J. P. Sch\"afermeyer, Berlin, Germany}
\email{jpschaefermeyer@gmail.com}

\author[D. S.]{David Sherry} \address{D. Sherry, Department of
Philosophy, Northern Arizona University, Flagstaff, AZ 86011, US}
\email{David.Sherry@nau.edu}

\numberwithin{equation}{section}

\begin{document}

\thispagestyle{empty}


\title[Klein \emph{vs} Mehrtens: restoring the reputation of a
modern]{Klein \emph{vs} Mehrtens: restoring the reputation of a great
modern}

\begin{abstract}
Historian Herbert Mehrtens sought to portray the history of
turn-of-the-century mathematics as a struggle of modern \emph{vs}
countermodern, led respectively by David Hilbert and Felix Klein.
Some of Mehrtens' conclusions have been picked up by both historians
(Jeremy Gray) and mathematicians (Frank Quinn).

We argue that Klein and Hilbert, both at G\"ottingen, were not
adversaries but rather modernist allies in a bid to broaden the scope
of mathematics beyond a narrow focus on arithmetized analysis as
practiced by the Berlin school.

Klein's G\"ottingen lecture and other texts shed light on Klein's
modernism.  Hilbert's views on intuition are closer to Klein's views
than Mehrtens is willing to allow.  Klein and Hilbert were equally
interested in the axiomatisation of physics.  Among Klein's credits is
helping launch the career of Abraham Fraenkel, and advancing the
careers of Sophus Lie, Emmy Noether, and Ernst Zermelo, all four
surely of impeccable modernist credentials.

Mehrtens' unsourced claim that Hilbert was interested in production
rather than meaning appears to stem from Mehrtens' marxist leanings.
Mehrtens' claim that [the future \emph{SS-Brigadef\"uhrer}] ``Theodor
Vahlen \ldots cited Klein's racist distinctions within mathematics,
and sharpened them into open antisemitism'' fabricates a spurious
continuity between the two figures mentioned and is thus an odious
misrepresentation of Klein's position.
\end{abstract}

\keywords{arithmetized analysis; axiomatisation of geometry;
axiomatisation of physics; formalism; intuition; mathematical realism;
modernism; Felix Klein; David Hilbert; Karl Weierstrass}

\maketitle
\tableofcontents

\epigraph{This may be regarded as a continuation of the \emph{Klein
Erlanger Programm}, in the sense that a geometrical space with its
group of transformations is generalized to a category with its algebra
of mappings.  Eilenberg--MacLane \cite[p.\;237]{EM} in 1945}

\section{Felix Klein}
\label{s1}

Historian Herbert Mehrtens sought to portray the history of
turn-of-the-century mathematics as a struggle of modern \emph{vs}
countermodern, represented respectively by David Hilbert and Felix
Klein; see e.g., Mehrtens \cite{Me90}.  Some of Mehrtens' conclusions
have been picked up by both historians (e.g., Jeremy Gray \cite{Gr08})
and mathematicians (e.g., Frank Quinn \cite{Qu12}).

To be sure, notable differences in outlook existed between Hilbert and
Klein.  Thus, Hilbert did not share Klein's intense interest in the
history of mathematics (see e.g., Rowe \cite[p.\;192]{Ro94}), and
Klein was less interested in axiomatics than Hilbert (see e.g., Weyl
\cite[p.\;16]{We85}).  However, such differences though undeniable
were extrapolated by scholars like Mehrtens to extravagant
proportions.

Against Mehrtens, we argue that Klein and Hilbert were not adversaries
but rather modernist allies in a bid to broaden the scope of
mathematics beyond a narrow focus on arithmetized analysis as
practiced by the Berlin school, so as to include set theory,
axiomatisation of physics, and other innovative directions.  To this
end, we analyze Klein's G\"ottingen lecture \cite{Kl96} and other
texts.

Felix Klein's 1895 G\"ottingen address was as influential as it is
apparently controversial, as we will see in later sections.

\subsection{The G\"ottingen address}  
\label{s31}

Here Klein spoke of
\begin{quote}
an important mathematical tendency which has as its chief exponent
Weierstrass\ldots{} I refer to the \emph{arithmetizing} of
mathematics.  \cite[p.\;241]{Kl96} (emphasis in the original)
\end{quote}
The passage makes clear Klein's appreciation of the significance of
the framework developed by Weierstrass and others.  According to
Klein, the framework constituted an advance over earlier reliance on
spatial intuition as a basis for proofs:
\begin{quote}
Gauss, taking for granted the continuity of space, unhesitatingly used
space intuition as a basis for his proofs; but closer investigation
showed not only that many special points still needed proof, but also
that space intuition had led to the too hasty assumption of the
generality of certain theorems which are by no means general. Hence
arose the demand for exclusively arithmetical methods of proof\ldots{}
(ibid.)
\end{quote}
The break with Gauss' view is particularly significant and underscores
the modernity of Klein's.  Such ``arithmetical methods'' meant that
\begin{quote}
nothing shall be accepted as a part of the science unless its rigorous
truth%
\footnote{The phrase ``rigorous truth' is meaningless according to
contemporary usage in modern mathematics, and is a mistranslation
found in Isabel Maddison's (on the whole adequate) translation.  Felix
Klein wrote, in remarkably modern phrasing: ``\ldots die Forderung
\emph{ausschlie{\ss}lich arithmetischer Beweisf{\"u}hrung [ist:] Als
Besitzstand der Wissenschaft soll nur angesehen werden, was durch
Anwendung der gew{\"o}hnlichen Rechnungsoperationen als
\emph{identisch richtig} klar erwiesen werden kann.}'' [emphasis
added] A correct translation of Klein's words is `\ldots the demand of
\emph{exclusively arithmetical proofs} [is:] only those propositions
are to be considered the secure possession of science which can
clearly be demonstrated to be identically valid by applying the usual
arithmetical operations.'  In particular, Klein's `identisch richtig'
(which Maddison rendered as `rigorous truth') is the German equivalent
of the English `identically valid'.}
can be clearly demonstrated by the ordinary operations of
analysis. (ibid.)
\end{quote}
Although vague notions, like magnitude, continuous variable, etc.,
were still in use in these new developments, Klein believed that
further refinements could be introduced through the limitations on the
notion of quantity (as in Kronecker's approach) or through the
application of symbolic language (the approach of Peano and his
school).  Klein continued:
\begin{quote}
Thus, as you see, while voluntarily acknowledging the exceptional
importance of the tendency, I do not grant that the arithmetized
science is the essence of mathematics\ldots{} I consider that the
essential point is not the mere putting of the argument into the
arithmetical form, but the mere rigid logic obtained by means of this
form.  (ibid., p.\;242)
\end{quote}
Klein felt that the pursuit of abstraction is an open-ended process
that need not stop with Weierstrass.  Klein's enthusiasm for the
arithmetization of analysis was evident.  Also in evidence is his
appreciation of new logical and (as we will show below) foundational
studies, even though the term \emph{logic} did not have the meaning we
attach to it today.%
\footnote{In this area Klein was more of an \emph{enabler} of new
mathematics than a direct \emph{contributor}; he expressed his
personal preference as follows: ``symbolic methods\ldots{} this
subject does not appeal to me personally'' (ibid., p.\;243).}
It is in this context that we should view Klein's further claim to the
effect that
\begin{quote}
it is not possible to treat mathematics exhaustively by the method of
logical deduction alone, but that, even at the present time,
\emph{intuition} has its special province. (ibid., p.\;242) (emphasis
added)
\end{quote}
We will deal with Klein's stance on intuition (mentioned in the
comment just cited) in Section~\ref{s22} immediately following.

\subsection{Spatial intuition}
\label{s22}

In the remainder of his address, Klein seeks to place the role of
spatial intuition in relation to logical and axiomatic developments.
Klein is somewhat ambiguous as to the meaning he attaches to the term
\emph{intuition}.  One can single out three possible meanings:
\begin{enumerate}
\item
intuition as an indispensable tool in research;%
\footnote{As when Klein writes: ``I might now introduce a historical
excursus, showing that in the development of most of the branches of
our science, intuition was the starting point, while logical treatment
followed'' \cite[p.\;246]{Kl96}.}
\item
intuition as an indispensable tool in teaching;%
\footnote{As when Klein writes: ``Among the teachers in our Gymnasia
the need of mathematical instruction based on intuitive methods has
now been so strongly and universally emphasized that one is compelled
to enter a protest, and vigorously insist on the necessity for strict
logical treatment.  This is the central thought of a small pamphlet on
elementary geometrical problems which I published last summer.  Among
the university professors of our subject exactly the reverse is the
case; intuition is frequently not only undervalued, but as much as
possible ignored'' \cite[p.\;248]{Kl96}.}
\item
\label{i2}
axiomatic accounts are insufficient and intuition must play its role
in mathematical arguments.%
\footnote{As when Klein writes: ``On the other hand I have to point
out most emphatically--and this is the negative part of my task--that
it is not possible to treat mathematics exhaustively by the method of
logical deduction alone, but that, even at the present time, intuition
has its special province'' \cite[p.\;242]{Kl96}.}
\end{enumerate}
What we wish to emphasize is that even if meaning \eqref{i2} occurs
here at all, it occurs on a sophisticated level as indicated by
Klein's endorsement of the arithmetic foundations for analysis
(therefore no more intuitive talk of real numbers) and his comments on
Green's theorem and electricity (emphasizing that physical intuition
is insufficient to prove mathematical theorems); see
Section~\ref{s23b}.  Klein's comment on Gauss quoted in
Section~\ref{s31} indicates that Klein clearly distanced himself from
a reliance on spatial intuition as replacement for analysis.  Even if
at some sophisticated level Klein thought that axiomatic approach will
be insufficient, no investigation of such a possible level of Klein's
term \emph{intuition} appears in Mehrtens's book, which contains
little indication that he would actually have the mathematical
wherewithal to carry out such an investigation.  The level of
mathematical competence possessed by Mehrtens is illustrated by his
comment that ``if it were possible to represent $\pi$ by an algebraic
equation, then the construction [i.e., squaring the circle] would be
possible with compass and straightedge'' \cite[p.\;111]{Me90}.  This
amounts to an incorrect claim to the effect that every algebraic
number is constructible ($\sqrt[3]{2}$ is algebraic but
non-constructible).

Klein felt that arithmetization of geometry meant that geometrical
objects are given by formulas and are dealt with by means of analytic
geometry and analytic methods:
\begin{quote}
The arithmetizing of mathematics began originally, as I pointed out,
by ousting space intuition; the first problem that confronts us as we
turn to geometry is therefore that of reconciling the results obtained
by arithmetical methods with our conception of space. By this I mean
that we accept the ordinary principles of analytical geometry, and try
to find from these the geometrical interpretation of the more modern
analytical developments.  \cite[p.\;243]{Kl96}
\end{quote}
The obvious advantage consists in the refinement of space intuition:
\begin{quote}
The net result is, on the one hand, a refinement of the process of
space intuition; and on the other, an advantage due to the clearer
view that is hereby obtained of the analytical results considered,
with the consequent elimination of the paradoxical character that is
otherwise apt to attach itself to them. (ibid., p.\;243)
\end{quote}
The next advantage is \emph{idealisation}, namely a mathematical form
of an imprecise intuition:
\begin{quote}
We ultimately perceive that space intuition is an inexact conception,
and that in order that we may subject it to mathematical treatment, we
\emph{idealize} it by means of the so-called axioms, which actually
serve as postulates. (ibid., p.\;243--244; emphasis added)
\end{quote}

\subsection{Idealisation in physics}
\label{s13}

As for idealisation in mechanics and mathematical physics, Klein
writes:
\begin{quote}
Throughout applied mathematics, as in the case of space intuition, we
must idealize natural objects before we can use them for purposes of
mathematical argument; but we find continually that in one and the
same subject we may idealize objects in different ways, according to
the purpose that we have in view. (ibid., p.\;244)
\end{quote}
In the realm of practical physics, idealisation provides precisely
defined objects like Green's function (see Section~\ref{s23b}), that
enable new physical insights as well as abstract mathematical
arguments.  Yet Klein's intuition is twofold:
\begin{quote}
\ldots{}[1] the cultivated intuition just discussed, which has been
developed under the influence of logical deduction and might almost be
called a form of memory; but rather of [2] the na\"ive intuition,
largely a natural gift, which is unconsciously increased by minute
study of one branch or other of the science.  (ibid., p.\;245--246)
\end{quote}
Klein elaborates as follows:
\begin{quote}
The word intuition (\emph{Anschauung}) is perhaps not well chosen; I
mean it to include that instinctive feeling for the proportion of the
moving parts with which the engineer criticises the distribution of
power in any piece of mechanism he has constructed; and even that
indefinite conviction the practiced calculator possesses as to the
convergence of any infinite process that lies before him.  I maintain
that mathematical intuition - so understood - is always far in advance
of logical reasoning and covers a wider field. (ibid.)
\end{quote}
Thus Klein's intuition [1] is \emph{cultivated} (under the influence
of logical deduction), while intuition [2] is the most basic in his
vision of mathematics, namely prelogical:
\begin{quote}
Logical investigation is not in place until intuition has completed
the task of idealisation. (ibid., p.\;247)
\end{quote}
Klein finds confirmation of the idea that intuition was the starting
point, whereas logical treatment followed, not only in the historical
origins of infinitesimal calculus, but also in Minkowski's development
of the theory of numbers.

As for the intuition [2] and pedagogy, Klein writes that
\begin{quote}
two classes at least of mathematical lectures must be based on
intuition; the elementary lectures which actually introduce the
beginner to higher mathematics - for the scholar must naturally follow
the same course of development on a smaller scale, that the science
itself has taken on a larger - and the lectures which are intended for
those whose work is largely done by intuitive methods, namely, natural
scientists and engineers. (ibid., p.\;248)
\end{quote}

To be sure, Klein's speculations on intuition are not profound; modern
philosophers may find them problematic.  Our goal here is not to argue
that Klein was a great philosopher but rather to indicate that his
preoccupations were those of modern mathematicians.  If Mehrtens
wishes to champion the cause of specifically \emph{modern}
mathematics, he cannot easily dismiss today's mathematicians as being
just as countermodern with regard to intuition as he claims Klein is.
Today mathematicians are still struggling with the role of intuition
in the creative process.

\subsection{Attitude toward logic and foundations}
\label{s23}

Klein recognized the significance of the contemporary developments in
logic associated with the names of Peano and others, when he spoke of
\begin{quote}
efforts made to introduce symbols for the different logical processes,%
\footnote{In the original: ``logischen Verkn{\"u}pfung''.  A better
translation would be ``logical connectives''.  While
\emph{Verkn{\"u}pfung} is often translated as \emph{operation} in
\emph{algebraic} contexts, the most accurate translation in
contemporary discussions of logic is \emph{connective}.}
in order to get rid of the association of ideas, and the lack of
accuracy which creeps in unnoticed, and therefore not allowed for,
when ordinary language is used.  In this connection special mention
must be made of an Italian mathematician, Peano, of Turin\ldots{}
\cite[p.\;242]{Kl96}
\end{quote}
Klein's warm relationship with Pasch (see Schlimm \cite{Sc13}) further
attests to Klein's visionary appreciation of contemporary developments
in axiomatic foundations.  At the same time, Klein voiced a cautionary
note: ``while voluntarily acknowledging the exceptional importance of
the tendency, I do not grant that the arithmetized science is the
essence of mathematics'' (see Section~\ref{s31} for a longer
quotation).  To Klein, the significance of the new methodology goes
hand in hand with a rejection of 19th century methodology:
\begin{quote}
From this outline of the new geometrical programme you see that it
differs greatly from any that was accepted during the first half of
this [i.e., 19th] century\ldots{} (ibid., p.\;244)
\end{quote}
Contrary to earlier ideas of mathematics as being a representation of
reality, Klein is clearly aware of the idealizing nature of
mathematical treatments, as shown by the passage quoted at the
beginning of Section~\ref{s13}.  As an example, Klein gives the
possibility of idealizing matter by either continuous or discrete
representations:
\begin{quote}
we treat matter either as continuous throughout space, or as made up
of separate molecules, which we may consider to be either at rest or
in rapid motion. (ibid., p.\;245)
\end{quote}
If matter could admit \emph{distinct} mathematical representations
according to Klein, then clearly Klein did not share the naive earlier
view of mathematics as a straightforward representation of reality.

\subsection{Klein on Green's function}
\label{s23b}

Klein goes on to note the distinction between, on the one hand, a
physical phenomenon and, on the other, a mathematical theorem
attempting to capture the latter.  He illustrates such a distinction
by citing the example that ``in electricity\ldots{} a conductor under
the influence of a charged point is in a state of electrical
equilibrium'' whose mathematical counterpart is the existence of
Green's function, and concludes:
\begin{quote}
You see here what is the precise object of these renewed
investigations; not any new physical insight, but abstract
\emph{mathematical argument} in itself, on account of the clearness
and precision which will thereby be added to our view of the
experimental facts.  (ibid.; emphasis added)
\end{quote}
It emerges that according to Klein the ultimate criterion of the
validity of a theorem is ``mathematical argument'' rather than
physical insight.

\subsection{Pedagogy}
\label{s16}

Going on to pedagogy, Klein notes: ``I must add a few words on
mathematics from the point of view of pedagogy'' (ibid., p.\;247).
University professors bear the brunt of Klein's critique:
\begin{quote}
Among the university professors of our subject exactly the reverse is
the case; intuition is frequently not only undervalued, but as much as
possible ignored. This is doubtless a consequence of the intrinsic
importance of the arithmetizing tendency in modern mathematics.  But
the result reaches far beyond the mark.  It is high time to assert
openly once for all that this implies, not only a false pedagogy, but
also a \emph{distorted view of the science}. (ibid., p.\;248)
(emphasis added)
\end{quote}
While cautious about possible deleterious effects on pedagogy, Klein
welcomed the arithmetizing tendency where appropriate, and made common
cause with Hilbert in most cases.

The visionary nature of Klein's \emph{Erlangen program} (EP) is widely
known and acknowledged even by Gray (see Section~\ref{s26}) though not
by Mehrtens (see Section~\ref{s2}).  A major weakness of Mehrtens'
book is his failure to address the importance of category theory in
modern mathematics.  The EP focused on transformations of objects
rather than objects themselves, a viewpoint recognizable as a
foundation rock of the category-theoretic approach; see Marquis
\cite{Ma09}.  What is remarkable is that even Mehrtens' nemesis
Bieberbach acknowledged that the EP was a precursor of the axiomatic
method (see Segal \cite[p.\;347, note\;56]{Se03}).  The EP furnishes
clear evidence in favor of Klein's modernism and lack of recognition
of this by Mehrtens constitutes massaging of evidence.

\subsection{Klein on infinitesimal analysis}

Klein's foresight of the eventual success of infinitesimal analysis in
modern mathematics is similarly remarkable.  Thus, Klein formulated a
criterion for what it would take for a continuum incorporating
infinitesimals to furnish a successful framework, in terms of the
availability of a mean value theorem in the framework; see Klein
\cite[p.\;213]{Kl08}.  Klein's criterion was endorsed by Abraham
Fraenkel \cite[pp.\;116--117]{Fran}.  Such a continuum was eventually
developed by Fraenkel's student Abraham Robinson \cite{Ro61}; see
Kanovei et al.\;\cite{18i} for details.

The proposal of Vinsonhaler \cite{Vi16}, Katz--Polev \cite{17h}, and
others for teaching calculus with infinitesimals relies on their
intuitive appeal and would have likely met with Klein's approval.
Modern frameworks for infinitesimal analysis occasioned a
re-evalulation of its history; see e.g., Bair et al.\;\cite{17a},
Bascelli et al.\;\cite{14a}, \cite{18c}, B\l aszczyk et
al.\;\cite{17d}.

\subsection{Courant on the \emph{Erlangen Program}}

Courant wrote:

\begin{quote}
This so-called \emph{Erlangen Program}, entitled `Comparative
Considerations on recent geometrical research' has become perhaps the
most influential and most-read mathematical text of the last 60 years.

Since the end of the 18th century, geometry had extraordinarily
thrived in France and Germany.  Alongside old elementary geometry and
analytical geometry, a large number of geometrical considerations had
[been] developed, which stood side-by-side unmotivated and without
mutual connections, and in this tangle [of ideas] even an expert could
hardly orient themselves anymore.  Klein felt the need to bring a
uniform ordering principle into this chaos, and he has solved this
task for the whole of geometry in way which one currently cannot
imagine coming nearer to being complete. The magic wand, with which
here Klein created order, was the group concept.

It [the group concept] permits to conceive every class of geometrical
investigations (like Euclidean and Non-Euclidean geometry, projective
geometry, line- and sphere-geometry, Riemannian geometry and topology)
as invariant theory relative to a given group of geometrical
transformations.  The \emph{Erlangen Program} constitutes for geometry
a similarly forceful ordering-principle as the periodic system is for
chemistry.

To this day, no geometrical theory can be considered finished if it
cannot clearly assert its place within the framework of the
\emph{Erlangen Program}.  Klein lived to have the satisfaction, fifty
years later,
%
%
to be able to very substantially contribute to the clarification of
the mathematical foundations of relativity theory, simply by,
essentially, applying his old thoughts from the \emph{Erlangen
program} to the new questions.  (Courant \cite[p.\;200]{Co26})
\end{quote}
Courant's comments clearly indicate the modern nature of Klein's EP.

\subsection{Hermann Weyl on Riemann surfaces}

Hermann Weyl comments as follows on the significance of the work of
Riemann on what are now called Riemann surfaces, and its clarification
by Klein:

\begin{quote}
It has to be admitted, of course, that Riemann himself slightly
disguised the true relationship between [complex] functions and
Riemann surfaces, by the form of his presentation---perhaps for the
one reason only that he did not want to visit overly alien ideas upon
his contemporaries; he disguised this relationship in particular by
only speaking of those multi-sheeted Riemann-surfaces covering the
plane which have finitely many ramification points, which [by the way]
even today are those which one primarily thinks of when the topic of
Riemann surfaces is mentioned, and that \emph{Riemann did not use the
more general idea (which only later was developed to transparent
clarity by K l e i n)}, an idea whose characteristic property can be
described as follows: any connection with the complex plane, and
generally any connection with three-dimensional point-spaces has been
severed, and severed in principle.  And yet, without any possible
doubt, \emph{only with Klein's conception} does the basic thought of
Riemann come to life and is the natural simplicity, vitality and
effective force of these ideas made visible.  The present work is
based on these [i.e., Klein's] thoughts.  (Weyl \cite[p.\;V]{We13};
emphasis added)
\end{quote}
Klein's role in Weyl's project was instrumental:
\begin{quote}
The most important of the basic thoughts [in Weyl's book] are due to
the man to whom I was permitted to dedicate this book in sincere and
profound reference.  Geheimrat%
\footnote{Often translated as \emph{privy councillor}.}
Klein has insisted, despite being overburdened by other tasks, and
despite his failing health, to discuss the entire subject matter in
frequent meetings with me; I owe a great debt of gratitude to him for
his remarks which on several occasions have caused me to replace my
initial presentation with a more correct and suitable one.  (ibid.,
p.\;IX)
\end{quote}
Weyl makes it clear that Klein made great \emph{sacrifices of time and
energy} to support the first \emph{abstract} and \emph{axiomatic} (and
hence textbook-modern) presentation of the theory of Riemann surfaces.
This is a factual argument against Mehrtens' counter-Klein thesis
analyzed in Section~\ref{s2}.

\section{Karl Weierstrass (1815--1897)}
\label{s210}

The seminal contribution of Karl Weierstrass to the development of
modern analysis is well known and requires no special comment.  In
this section we will focus on Weierstrass' comments on intuition and
ethnicity, and compare them to Klein's.

\subsection{Letter to Kovalevskaya}
\label{s21}

Sanford Segal notes that Weierstrass' 1883 letter to his student
Kovalevskaya was cited by SA footsoldier Teichm\"uller in a 1938
lecture; see Segal \cite[p.\;405]{Se03}.  In his letter, Weierstrass
expressed himself as follows:
\begin{quote}
Among the older mathematicians there are different sorts of human
beings; this is a trivial proposition, which nevertheless explains
much.  My dear friend Kummer, for example, was not interested in what
happened in mathematics as a whole, neither while he applied all his
energy to find the proofs of the higher reciprocity laws, nor, and
then less than ever, after he had expended his energies in these
pursuits.  His attitude [to new developments in mathematics] is, if
not dismissive, indifferent.  If you tell him that Euclidean Geometry
is based on an unproved hypothesis, then he [readily] admits this;
[but] to proceed from this insight to asking the question `How does
geometry develop without this hypothesis?', now that is contrary to
his [mental] nature, and the exertions aimed at this question, and the
general investigations which follow from this question and which
liberate themselves from what is empirically given or assumed, are
moot speculations to him, or even an abomination.  (from an 1883
letter from Weierstrass to Kovalevskaya reproduced in Mittag-Leffler
\cite[p.\;148--150]{Mi02}; translation ours)
\end{quote}
Having thus analyzed Kummer, Weierstrass goes on to provide his
insights into Kronecker's temperament:
\begin{quote}
Kronecker is different, he familiarizes himself quickly with
everything that is new, his facility of perception enables him to do
so, but it does not happen deeply; he does not have the talent to
apply the same interest to good work of others as to his own work.
(ibid.)
\end{quote}
Weierstrass then proposes the following analysis of such differences
in mathematical temperaments:
\begin{quote}
This is compounded by a defect which can be found in many very
intelligent people, especially those from the \emph{semitic tribe}[;]%
\footnote{Rowe translates \emph{Stamm} as ``stock'' but ``tribe''
seems both more accurate and more consistent with Weierstrass' tone in
the letter.}
he [i.e., Kronecker] does not have sufficient imagination (I should
rather say: \emph{intuition}) and it is true to say that a
mathematician who is not a little bit of a poet, will never be a
consummate mathematician.  Comparisons are instructive: the
all-encompassing view which is directed towards the Highest, the
Ideal,%
\footnote{Rowe translates Weierstrass' phrase ``allumfassende auf das
h{\"o}chste, das Ideale gerichteten Blick'' as ``all-embracing vision
focused on the loftiest of ideals'' but the translation is incorrect.
Weierstrass is \emph{not} using `h{\"o}chste' as an \emph{adjective}
modifying `Ideale'; rather, Weierstrass uses a juxtaposition of a
\emph{neuter nominalized superlative} (`das h{\"o}chste') with a
\emph{neuter nominalized adjective} (`das Ideale'), separating the two
by a comma instead of the conjunction `and', to create a rhetorical
effect in his letter.}
marks out, in a very striking manner, Abel as better than Jacobi,
marks out Riemann as better than all his contemporaries (Eisenstein,
Rosenhain), and marks out Helmholtz as better than Kirchhoff (even
though the latter did not have a droplet of \emph{semitic blood}).%
\footnote{In the original: ``Dazu kommt ein Mangel, der sich bei
vielen h{\"o}chst verst{\"a}ndigen Menschen, namentlich bei denen
semitischen Stammes findet[;] er besitzt nicht ausreichend Phantasie
(Intuition m{\"o}chte ich lieber sagen) und es ist wahr, ein
Mathematiker, der nicht etwas Poet ist, wird nimmer ein vollkommener
Mathematiker sein. Vergleiche sind lehrreich: Der allumfassende auf
das h{\"o}chste, das Ideale gerichtete Blick zeichnet Abel vor Jacobi,
Riemann vor allen seinen Zeitgenossen (Eisenstein, Rosenhain),
Helmholtz vor Kirchhoff (obwohl bei dem letztern kein Tr{\"o}pfchen
semitischen Blutes vorhanden) in ganz eclatanter Weise aus.''}
(ibid.; emphasis added)
\end{quote}
Notes Rowe: ``These remarks, echoed in the very same stereotypes set
forth by anti-Semites like Dingler and Bieberbach fifty years later,
are a good illustration of how deeply rooted such thinking was in
German culture'' (Rowe \cite[p.\;443]{Ro86}).

Weierstrass' speculations as to the insufficient endowment in
\emph{intuition} on the part of the \emph{semitic tribe} marks out
Weierstrass as a better candidate than Klein (who never proffered such
speculations) for countermodern leadership, by Mehrtens' own
standards.  We would like to suggest that Mehrtens' book is critically
flawed in having misidentified its protagonists: while Hilbert can
serviceably provide modern leadership, his opposite numbers at the
helm of the \emph{intuitive countermoderns} should have been
Weierstrass and his students (see Bair et al.\;\cite{18b}) rather than
Klein.

We would like to clarify that Weierstrass made these comments in a
private letter, and did not anticipate that his classification will be
amplified through publication by Mittag-Leffler.  Klein, on the other
hand, made his comments on Teutonic, Latin and Hebrew characteristics
(see Section~\ref{s46}) in a public address (Klein \cite{Kl94}).  It
is unclear whether Mehrtens would characterize Weierstrass' comments
as \emph{biologistisch-rassistische}.%
\footnote{See note~\ref{f12} in Section~\ref{s29}.}

\subsection{Weierstrass and Klein on Kronecker}
\label{s22b}

In Section~\ref{s21} we observed that Karl Weierstrass attributed
perceived shortcomings of Kronecker's mathematical outlook to the
latter's membership in the \emph{semitic tribe} with an attendant
deficiency in \emph{intuition}.  Weierstrass' remarks on Kronecker can
be profitably compared with those made by Klein, in a letter to
Friedrich Althoff, ``the kingpin of the Prussian university system''
(see Rowe \cite[p.\;424]{Ro86}), concerning mathematics at Berlin:
\begin{quote}
Without question the positive aspects have been borne primarily by
Kronecker.  In this respect I must not withhold my praise.  That
Kronecker, even his last years of life, was able to bring new ideas to
our science with such youthful ambition, and thereby to uphold
Berlin's old fame as a center for mathematical research in a new form,
that is an accomplishment one can only admire without reservation.  My
critique merely concerns the one-sidedness with which Kronecker, from
a philosophical standpoint, fought against various scientific
directions that were remote from his own\ldots{} This one-sidedness
was probably less grounded in Kronecker's original talents than it was
in the disposition of his character.  Unconditional mastery, if
possible over all of German mathematics, became more and more the goal
which he pursued with all the cleverness and tenacity he could muster.
Little wonder that there is no one to take his place now that he has
left the arena.  (Klein to Althoff, as quoted in Rowe
\cite[p.\;442]{Ro86})
\end{quote}
Regardless of whether one agrees with Klein's assessment of
Kronecker's mathematics and/or politics, one has to acknowledge that
from the viewpoint of racial harmony, Klein's remarks are
unobjectionable; Kronecker's relation to what Weierstrass referred to
as the \emph{semitic tribe} is not mentioned.  Whenever Klein did
mention a mathematician's Jewishness, it was in laudatory terms like
the following:
\begin{quote}
Personally, [James Joseph] Sylvester was extremely engaging, witty and
effervescent.  He was a brilliant orator and often distinguished
himself by his pithy, agile poetic skill, to the mirth of everyone. By
his brilliance and agility of mind he was a genuine representative of
his race; he hailed from a purely Jewish family, which, having been
nameless before, had adopted the [sur]name Sylvester only in his
generation.  (Klein \cite[p.\;163]{Kl26})
\end{quote}
Klein's philosemitic comments in this area stand in contrast with
those penned by Weierstrass and cited in Section~\ref{s21}, as well as
those penned by Vahlen and cited in Section~\ref{s310}.

\subsection{Letter to Schwarz}
\label{1234}

Weierstrass opposed the looming appointment of Sophus Lie at Leipzig
on the grounds that the foreigner Lie was not of such consequence as
to warrant a passing over of actual \emph{countrymen} (i.e.,
compatriots); see Stubhaug \cite[p.\;317]{St02}.  Weierstrass went on
to express his scorn in a 1885 letter written from Lake Geneva to his
student Schwarz at G{\"o}ttingen, in the following terms:
\begin{quote}
Du Bois-Reymond really sometimes nails it; years ago already he called
the \emph{Trifolium} of Klein--Lie--Mayer an `acolyte society.%
\footnote{In the original: ``Du Bois-Reymond trifft doch zuweilen den
Nagel auf den Kopf, er nannte vor Jahren schon das Trifolium
Klein--Lie--Mayer `soci{\'e}t{\'e} thurif{\'e}raire'.''}
(Weierstrass as cited in Confalonieri \cite[p.\;288]{Co13})
\end{quote}
The \emph{nailing} comment is followed by a piece of advice for
Schwarz:
\begin{quote}
But, dear friend, don't let yourself be led astray, go on to correctly
walk your walk, teach and work in your thorough style---thus you will
best counteract the swindle,%
\footnote{\label{f12b}In the original: ``Schwindel''.  Weierstrass'
term apparently refers to the activity of the Klein--Lie--Mayer
society and connects well with Weierstrass' epithet \emph{Blender}
used by Weierstrass in 1892 in reference to Klein.}
even though \emph{initially}, you will have to brace yourself for
this, the big crowd%
\footnote{The words ``der gro\ss e Heifer'' in Confalonieri's
transcription are possibly a corruption of ``der gro\ss e Haufen''.
Weierstrass' manuscript shows a black speck above the squiggly line
representing this word.  The speck is most likely an overline that was
once used in German to distinguish a handwritten `u' from a
handwritten `n'.  Confalonieri seems to have mistaken the speck for a
dotted `i', leading him to guess `Heifer' (not a word in German).}
will turn to the newly rising sun.%
\footnote{The words ``der neu mitgehenden Sonne'' in Confalonieri's
transcription are possibly a corruption of ``der neu aufgehenden
Sonne''.  Note that the word in question starts with the glyph that
Weierstrass uses for `a'.}
(Weierstrass as cited in Confalonieri \cite[p.\;288]{Co13})
\end{quote}
Weierstrass appears to be mixing metaphors here.  The sunflower
follows the sun, just as a crowd blindly follows a new fad.
Weierstrass' \emph{newly rising sun} appears to refer to Klein.

These comments by Weierstrass are significant because they indicate
that he had a fundamental disagreement with Klein over the nature and
future of mathematics.  Weierstrass seems to have viewed Klein's
interests such as the \emph{Erlangen program} as a passing fad (that
goes away as quickly as the sun sets).  Huygens had similar feelings
about Leibniz' calculus, in that he thought it was merely a
repackaging of existing techniques.

\section{Mehrtens on Klein}
\label{s2}

In november 1933 philosopher Hugo Dingler filed a pamphlet against
Klein (1849--1925). 

\subsection{Dingler's pamphlet; J-type and S-type}
\label{s31b}

In his pamphlet, Dingler claims that Klein was half-Jewish, that Klein
filled G\"ottingen with Jews and foreigners, and that Klein was hungry
for power to control German mathematics and re-make it along Jewish
lines and in sum was \emph{un-German}; see Rowe\;\cite{Ro86}.

\subsection{Bieberbach versus Dingler}
\label{unsavory}

A quarter century after defending his thesis under Klein, Ludwig
Bieberbach made a claim \emph{contrary} to Dingler's regarding Klein,
with the latter now becoming emphatically \emph{German}; see
Bieberbach \cite{Bi34}.  Thus Bieberbach sought to attribute his own
(Bieberbach's) views concerning \emph{German mathematics} to Klein
himself, and to co-opt the latter in the service of an unsavory
ideology.

Bieberbach relied on a dichotomy of \emph{S-type} versus \emph{J-type}
borrowed from Jaensch \cite{Ja31}.  Here \emph{S-type} (for
\emph{Strahltypus}) refers to a ``radiating'' type that ``only values
those things in Reality which his intellect infers in it'' and
moreover ``denies the connection to an outer reality that is not
mentally constructed'' (see Segal \cite[p.\;365]{Se03}).

In contrast, the \emph{J-type} (or \emph{I-type}, for
\emph{Integrationstypus}) refers to one who ``is wide-open to
Reality'' and ``lets the influence of experience stream into him''
(see Segal \cite[362--363]{Se03} for details).

The Jews Jacobi and Landau among others was construed as
\emph{S-types}, whereas Klein as a \emph{J-type}.  On occasion,
Bieberbach used the terms \emph{Aryan} and \emph{non-Aryan}
\cite[p.\;177]{Bi35}; for further discussion see Segal
\cite[p.\;380]{Se03}.

\subsection{Mehrtens' choice; C-type and M-type}
\label{s33}

Over half a century later, historian Herbert Mehrtens chose to back
Bieberbach in the Dingler--Bieberbach disagreement over Klein.%
\footnote{This is not to imply that Mehrtens endorsed either Dingler's
or Bieberbach's views in the sense of their unsavory political or
philosophical orientations.  Nor is it an assertion that Mehrtens' own
political or philosophical orientations are \emph{not} unsavory.}
The historian confidently announced that ``Klein is indeed a
representative of counter-modernism'' (Mehrtens \cite[p.\;520]{Me96}).
What Mehrtens found countermodern in Klein is Klein's mathematical
\emph{outlook} rather than his mathematical \emph{output}.

To contrast Klein with Hilbert, Mehrtens introduced a dichotomy of
\emph{countermodern} (\emph{C-type} for short) and \emph{modern}
(\emph{M-type} for short), with Klein construed as \emph{C-type} and
Hilbert as \emph{M-type} (\emph{C-} and \emph{M-} notation ours).

Mehrtens not only claimed that Klein was \emph{C-type} but also
dropped both damaging innuendo that Klein may have ultimately been an
enabler of the national-socialist ideology and odious allegations of
racism, as analyzed in Sections~\ref{s32b} and \ref{s29}.

Note that the onus is on Mehrtens to convince the reader that he is
pursuing a meaningful dichotomy of \emph{C-type} versus \emph{M-type}.
We argue that by pursuing extraneous issues, Mehrtens obscures the
true vital issues of the period, and moreover reveals his ideological
\emph{parti pris} in the process.

\subsection{Racism}
\label{s32b}

Mehrtens' page~520 contains no fewer than six occurrences of the term
\emph{racism}; thus, we find:
\begin{quote}
The case of Bieberbach as well as the behaviour of some other
mathematicians suggest the thesis that mathematical counter-modernism
is correlated to nationalism and eventually also to \emph{racism}.
(ibid.; emphasis added)
\end{quote}
Furthermore,
\begin{quote}
Klein is indeed a representative of counter-modernism, but his
\emph{racism} is more of the theoretical type.  (ibid.; emphasis
added)
\end{quote}
Mehrtens' procedure here is objectionable on several counts:
\begin{itemize}
\item
The juxtaposition of the two passages (one on Bieberbach and one on
Klein) in close proximity suggests a spurious affinity between them;
\item
the claim of being a ``racist of a theoretical type'' presupposes
being a racist in the first place, an allegation against Klein that
Mehrtens has not yet established.
\end{itemize}
Now the term \emph{racist} can be given (at least) two distinct
meanings:
\begin{enumerate}
\item someone interested in analyzing differences in intellectual
outlook among distinct ethnicities (\emph{racist1});
\item someone who believes in the inferiority of one ethnicity to
another based on such differences, and advocates corrective action
(\emph{racist2}).
\end{enumerate}
Klein's discussions of ethnic differences possibly make him
\emph{racist1}.  Yet Mehrtens' comments on Bieberbach clearly indicate
that Mehrtens has \emph{racist2} in mind when he uses the term:
\begin{quote}
When the National Socialists came to power in 1933, [Bieberbach]
attempted to find political backing for his counter-modernist
perspective on mathematics, and declared both, intuition and
concreteness, to be the inborn characteristic of the mathematician of
the German race, while the tendency towards abstractness and
unconcrete logical subtleties would be the style of Jews and of the
French (Mehrtens 1987).  He thus turned counter-modernism into
outright \emph{racism} and anti-modernism.  (Mehrtens
\cite[p.\;519]{Me96}; emphasis added)
\end{quote}
Thus the \emph{racism} case against Klein as found in Mehrtens
\cite{Me96} is based on equivocation on the meaning of the term
\emph{racist}.  The remaining four occurrences of the term
\emph{racism} on page 520 in Mehrtens are less tendentious; yet
Mehrtens' reader can well wonder why the issue is being discussed in
such detail at all, if pinning a \emph{latent racism} slur on an
allegedly \emph{C-type} Klein were not one of Mehrtens' intentions.
Mehrtens goes so far as to imply a connection between Klein's comments
and the \emph{Nationalsozialismus}; see Section~\ref{s29}.

The historical record indicates that Klein struggled valiantly to hire
Jewish mathematicians like Hurwitz and Schoenflies, and conducted a
warm correspondence with Pasch, Gordan, and others.  Klein at times
put his own reputation on the line to do so, as in his correspondence
with Althoff concerning the promotion of Schoenflies.%
\footnote{Schoenflies went on to do work fundamental for both group
theory and modern spectroscopy.  This work is widely used to this day.
His correspondence with Fedorov aimed at correcting errors in earlier
versions of their classification of crystallographic groups is one
famous early example of the modern phenomenon of mathematicians from
different linguistic communities collaborating incrementally to arrive
at a valid result; see Schwarzenberger \cite[pp.\;162--163]{Sc84}.}

Mehrtens' strategy in the face of such facts is 
\begin{enumerate}
\item
to acknowledge that Klein appreciated some Jewish mathematicians, but
\item
to claim that this does not contradict Mehrtens' thesis that Klein,
perhaps unwittingly and without foreseeing the criminal abuses that
were to come, produced texts and personal statements which were grist
to the mill of antisemitic tendencies in Germany at the time
\cite[p.\;217]{Me90}.
\end{enumerate}
This stance of Mehrtens' amounts to assignment of guilt by association
and therefore represents a character smear against Klein.

\subsection{Internationalism}

Mehrtens goes on to raise the issue of internationalism:
\begin{quote}
The converse proposition would be that modernism is related to
liberalism and internationalism. Indeed, the leading representative of
the modem style, David Hilbert, can be rated as a liberal and an
internationalist.  (Mehrtens \cite[p.\;520]{Me96})
\end{quote}
A reader would gather from this comment that Klein may perhaps
\emph{not} be rated as a liberal, being a \emph{C-type} contrary to
Hilbert's \emph{M-type} (see Section~\ref{s33}).  In his review of
Mehrtens' approach, Albert Lewis aptly remarked:
\begin{quote}
The two world wars and their aftermaths brought out the stark contrast
between \emph{international} and \emph{national} mathematical biases
while the center of gravity of mathematical activity moved towards the
United States. \; (Lewis\;\cite{Le02}; emphasis added)
\end{quote}
Indeed, relative to the \emph{national} {vs} \emph{international}
dichotomy mentioned by Lewis, Felix Klein's outlook clearly belongs to
the latter (as correctly pointed out by Dingler; see
Section~\ref{s31b}), which is the opposite of the countermodern role
Mehrtens seeks to pin on Klein.

\subsection{What is mathematical modernism? Mehrtens '93}

Herbert Mehrtens and other historians ``have attempted to give a
Marxist analysis of the connection between mathematics and productive
forces, and there have been philosophical studies about the
communication processes involved in the production of mathematical
knowledge'' (Mehrtens et al.\;\cite[pp.\;ix--x]{Me81}).  In a 1993
text, Mehrtens posits that
\begin{quote}
a scientific discipline exchanges its knowledge products plus
political loyalty in return for material resources plus social
legitimacy.  (Mehrtens \cite[p.\;220]{Me93})
\end{quote}
Readers concerned that that the tools of a crude marxism may not
provide sufficient discrimination to deal with delicate issues of the
transformation that took place in mathematics at the beginning of the
20th century, will have their apprehensions confirmed when they reach
Mehrtens' definition of modernism:
\begin{quote}
By 1930 it was quite clear what the term \emph{modern} meant when
applied to mathematics: the conceptual study of abstract mathematical
concepts characterized by axioms valid for sets or other undefined
elements and presented by proceeding from the elementary concepts to
the more complicated structures; a hierarchical system of mathematical
truths rigorously proved, the language applied having hardly any other
function than to label the objects and to ensure the validity of
statements.  ``Modern'' mathematics in this sense had no
extramathematical \emph{meaning}, did not indicate possible fields or
objects of application, was devoid of hints to the historical or
heuristical background of the theory, and at most was in a highly
implicit manner structured along didactical guidelines.  (Mehrtens
\cite[p.\;224]{Me93}; emphasis added)
\end{quote}
Among the problematic aspects of Mehrtens' definition of modernism are
the following:
\begin{enumerate}
\item
Mehrtens makes no distinction between the levels of language and
metalanguage, or mathematics and metamathematics, which are essential
to understanding Hilbert's finitist program;
%
%
\item
Mehrtens rules out applications by definition (in his second
sentence), as well as didactic concerns (little wonder Klein didn't
fare well in Mehrtens' book);
\item
\emph{meaning} is ruled out by definition, which would exclude Hilbert
from Mehrtens' modernism;
%
%
\item
Mehrtens' dismissive comment on history (``devoid of hints to the
historical \ldots{} background of the theory'') flies in the face of
the fact that Hilbert's \emph{Grundlagen der Geometrie} explicitly
refers to Euclid's \emph{Elements}.
\end{enumerate}
To elaborate on the last point, note that Hilbert deals not merely
with the axioms but also with Euclid's theory of area and theory of
proportion.  Hilbert incorporates the perspective of Euclid as an
organic part of his presentation in the \emph{Grundlagen} (see Hilbert
\cite{Hi99}).  Mehrtens goes on to claim that
\begin{quote}
The central cognitive definition of mathematics by its \emph{pure,
modern core} was related to those task-oriented fields [of applied
mathematics and mathematical pedagogy].  (Mehrtens
\cite[p.\;226]{Me93}; emphasis added)
\end{quote}
However, the positing of a ``pure, modern core'' by Mehrtens remains a
mere assumption echoing those made two pages earlier.  Mehrtens
continues with a thinly veiled dig against Klein's book \cite{Kl08}:
\begin{quote}
Lectures on ``elementary mathematics from higher standpoint'' served
\ldots{} to relate the pedagogical branch to the development of
scientific mathematics, and to attempt to axiomatize mathematical
theories for physics aimed at an immediate relation between this field
of application and pure mathematics.  (Mehrtens \cite[p.\;226]{Me93})
\end{quote}
Again, the exclusion of applications from the ``pure, modern'' realm
is \emph{assumed} rather than \emph{argued} by Mehrtens.  His reader
may be justifiably shocked to learn that
\begin{quote}
In a lecture in 1926 [Bieberbach] sharply attacked [David] Hilbert,
the dean of mathematical modernity, and depicted modern theories as
`skeletons in the sand of the desert of which nobody knows whence they
come and what they have served for.'  (Mehrtens \cite[p.\;227]{Me93})
\end{quote}
However, the suitably shocked reader may well feel it to be a weakness
of Mehrtens' depiction of Klein as a \emph{C-type}, since obviously no
similar anti-Hilbert lecture was ever given or imagined by Klein.
Mehrtens wants us to believe that
\begin{quote}
the social system of mathematics was (and is) \emph{interested not in
meaning but rather in production}.  The competing formalist style was
obviously productive, and the corresponding foundational research
program, Hilbert's `metamathematics,' imposed no restrictions on
existing mathematical theories except the demand for logical
coherence.  (Mehrtens \cite[p.\;230]{Me93}; emphasis added)
\end{quote}
Was Hilbert in fact interested in \emph{production} rather than
\emph{meaning} as Mehrtens claims here?  Mehrtens probably believes
that, too.  However, he presents no evidence to support such a
sweeping claim that is more revealing of his (Mehrtens') ideological
commitments than of Hilbert's views, as is Mehrtens' comment to the
effect that
\begin{quote}
In [Hilbert's] program, mathematics is of the highest generality, and
\emph{meaning} and utility of new knowledge are hardly visible
anymore.%
\footnote{In the original: ``In seinem Programm ist die Mathematik von
h{\"o}chster Allgemeinheit, und Sinn und Nutzen neuen Wissens sind
kaum mehr zu erkennen.''}
(Mehrtens \cite[p.\;380]{Me90}; emphasis added)
\end{quote}
On Hilbert's commitment to meaning and utility, \emph{pace} Mehrtens,
see further in Section~\ref{s315}.

\subsection{What is mathematical modernism? Mehrtens '96}
\label{s25}

In 1996 Mehrtens provides an additional definition of mathematical
modernism: ``The counter-modernist attitude arises with modernism''
\cite[p.\;522]{Me96} and
\begin{quote}
It is part of modernity of the modern world.  Turning to the political
side of my topic, I want to state very briefly the sociological
argument.  I take mathematical modernism to be the defining center of
a socially modernized professional and autonomous mathematics
business, concentrating and symbolizing itself as a research
discipline.  Modernism is, so to speak, the avant-garde of the
profession, defining the real, the pure, and the most progressive type
of mathematics.  The business, however, includes much more than
avant-garde research.  Mathematicians of a more conservative type or
\emph{mathematicians in applications}, hybrid fields or in teaching
cannot fully identify with the modernist attitude. (Mehrtens
\cite[p.\;523]{Me96}; emphasis added)
\end{quote}
The problem with such a definition is that Klein becomes
\emph{tautologously} a \emph{C-type} here, because ``mathematicians in
applications" are \emph{defined} as incapable of fully identifying
with the modernist attitude.  The talk about ``socially modernized
professional and autonomous mathematics business'' would indicate on
the contrary that Klein is a modernist, since nobody contributed to
the flourishing of professional mathematics institutions as he did.

Readers can wonder what ``the real, the pure, and the most progressive
type of mathematics'' is exactly.  Mehrtens defines modernist
mathematics as \emph{the real, the pure, and the most progressive
type} but the definition seems circular.  What is \emph{progressive}
mathematics if not \emph{modern}?

Readers can have similar doubts about the \emph{a priori} exclusion of
applied mathematicians from the rarefied realm of the modernist
species.  Thus, Abraham Robinson published numerous articles in
applied mathematics as well as a joint book called \emph{Wing theory}
\cite{Ro56} (yes, these are \emph{airplane} wings, not a branch of
\emph{ring theory}).  Yet his philosophical stance is close to
Hilbert's formalism; see Robinson's essay ``Formalism\;64''
\cite{Ro65} which would normally gain him admittance to the select
modernist club if Mehrtens' enthusiasm for David Hilbert's formalism
is any guide.

Even more problematic for Mehrtens is the fact that Hilbert, like
Klein, was specifically interested in axiomatizing physics:
\begin{quote}
Another author who influenced Hilbert deeply during this K\"onigsberg
period was Heinrich Hertz, as reflected in [Hilbert's] sixth Paris
problem, which pointed toward the axiomatization of mathematical
physics.  (Rowe \cite[p.\;180]{Ro17})
\end{quote}
(See also Corry \cite{Co04a}.)  If applications are specifically
excluded from the realm of mathematical modernism, what is one to make
of this application to physics pursued by Mehrtens' modernist hero
Hilbert?

When Mehrtens proclaims: ``I take mathematical modernism to be the
defining center of a socially modernized\ldots{} mathematics
business,'' anyone with a background in (mathematical) logic might
wonder how one can define mathematical \emph{modernism} in terms of
\emph{modernized} business without committing an elementary logical
error of circularity.  In short, what Mehrtens provides here is not a
definition of modernism but rather an assortment of cherry-picked
conditions designed specifically to exclude\ldots{} Felix Klein.  Such
a technique could be described in popular parlance as moving the
goalposts to score a point.

\subsection{Stereotyping}
\label{s26c}

Does it ever occur to Mehrtens that significant \emph{differences} may
exist among the alleged \emph{C-types}?  The answer is that it does,
as when he speaks of
\begin{quote}
the individual who constructs himself as part of a higher order, and
who, from that higher order, \ldots{} receives the gift of a
mathematical talent and thus the gift of insight into that higher
order.  This applies to Poincar\'e, Klein, Bieberbach and others, but
not so much to the more radical and pessimistic Brouwer.  \emph{I
shall not discuss the individual differences here.}''  (Mehrtens
\cite[pp.\;525--526]{Me96}; emphasis added).
\end{quote}
Mehrtens' passage here is assorted with a footnote\;6 offering
tantalizing glimpses into what such undiscussed differences might be:
\begin{quote}
6. But I would like to recall the `polythetic' character of concepts
like `fundamentalism'.  One could, tentatively, divide the family
resemblances into a `technocratic', progress-oriented group (Klein,
Poincar\'e), a traditionalist or `mandarin' group (Kronecker), and a
romantic, mystical variant, connecting mathematical fundamentalism
with a critique of progress and civilization as Brouwer did.
(Mehrtens \cite[p.\;526, note\;6]{Me96})
\end{quote}
Alas, the reader never learns whether Bieberbach was fundamentalist,
mandarin, mystical, or technocratic; perhaps all of the above.  We
therefore have the following classification into $C$-subtypes and
their leading representatives following Mehrtens:
\begin{enumerate}
\item[$C_1$] Bieberbach (fundamentalist?);
\item[$C_2$] Brouwer (mystical);
\item[$C_3$] Kronecker (mandarin);
\item[$C_4$] Klein, Poincar\'e (technocrat).
\end{enumerate}
We note that, \emph{pace} Mehrtens, stereotyping mathematicians into
such categories serves no useful purpose if the historian's goal is
\emph{meaningful} history, though it may well serve a goal of the
\emph{production} of marxist historiography.  In the end Mehrtens'
claim concerning Hilbert that he was interested in production rather
than meaning applies only to the author of the claim.

Mehrtens' systematic stereotyping reminds one of nothing more than
similar procedures adopted by his nemesis Bieberbach.  Thus, Segal
summarizes Bieberbach's racial theorizing as follows:
\begin{quote}
Gauss was also contrasted with Carl Gustav Jacobi (a Jew).  Jacobi was
``oriental" and had a ``heedless will to push through his own
personality." Gauss was characterized as ``nordisch-falisch," a term
borrowed from H.F.K. Gunther's racial theories; similarly, Euler was
``ostisch-dinarisch," another similar term.  (Segal
\cite[p.\;363]{Se03})
\end{quote}
See there for a discussion of Bieberbach's classification of
\emph{J-types} into subtypes $J_1, J_2, J_3$.  Mehrtens' final
paragraph contains the following gem:
\begin{quote}
When the ideology of sound, realist common sense had crossed the
border to Germany in the early thirties, mathematical modernism took
the opposite route with a group of scholars that were to name
themselves ``Bourbaki'' and to become the last \emph{high priests} of
mathematical modernism before the \emph{post-modern} era.  (Mehrtens
\cite[p.\;527]{Me96}; emphasis added)
\end{quote}
It is unclear what meaningful historiographic purpose is served by
stereotyping the Bourbaki as ``high priests.''  

Mehrtens' claim that in mathematics, \emph{modernism} was followed by
something called \emph{post-modernism} may strike many a reader as
novel.  Then again such a reader may be unfamiliar with the
intricacies of the received academic lingo.  What objective does it
serve to postulate a specifically \emph{post}-modern phase in the
development of mathematics?  We will venture an explanation in
Section~\ref{s26b}.

\subsection{\emph{Late capitalism} and Kramer's diagnosis}
\label{s26b}

Writes Gray:
\begin{quote}
Mehrtens's critique was written in a post-Marxist spirit, influenced
by such writers as Foucault.  Modernization, for him, is not progress,
it is also part of the catastrophe of Nazism.  If he is less clear
that the search for meaning and a place in the world has its good
side, it is only because he sees more clearly the ways in which
\emph{late capitalism} is antithetical to all of that.  (Gray
\cite[p.\;10]{Gr08}; emphasis added)
\end{quote}
Here Gray clearly acknowledges the marxist source of Mehrtens'
inspiration.  Gray goes on to state:
\begin{quote}
Two further avowedly speculative chapters close [Mehrtens'] book,
which go further into \emph{cultural criticism} than I need to follow
here.''  (ibid., footnote\;13; emphasis added)
\end{quote}
The root of the dilemma, as far as \emph{cultural criticism} is
concerned, is that Mehrtens and similar-minded marxist academics have
a problem with the symbiotic relationship between capitalism on the
one hand, and modernism and high culture on the other.  After its
revolutionary beginnings as a radical movement,%
\footnote{In this connection, Mehrtens writes that \emph{modernity}
can be a \emph{radical} stance: ``die m{\"o}gliche Radikalit{\"a}t der
Moderne'' \cite[p.\;182]{Me90}.}
modernism went on to enjoy a symbiotic relationship with what Gray
refers to as \emph{late capitalism}, providing the source of the
enmity toward modernism on the part of marxist academics disenchanted
with bourgeois society.  The said academics felt betrayed by
modernism, and can therefore speak approvingly only of
\emph{post}modernism (or pop art), never of modernism itself, as
poignantly summarized by Hilton Kramer in his essay `Modernism and its
enemies':
\begin{quote}
\ldots{} This view may be summed up as follows: Modernism claimed to
be revolutionary, it claimed to be anti-bourgeois, it promised us a
new world, but it turned out to be a coefficient of bourgeois
capitalist culture, after all, and we therefore reject the claim of
high culture and must work to destroy the privileged status it enjoys
in the cultural life of the bourgeois democracies.  (Kramer
\cite[p.\;13]{Kr88})
\end{quote}
To Mehrtens, \emph{applied} mathematics is similarly a ``coefficient
of bourgeois capitalist culture'' aiding and abetting not merely the
bourgeois democracies but the NS regime as well:
\begin{quote} 
The plan for an international congress for mechanics in Germany could
be sold to the aircraft ministry as a necessity for productive
aircraft research.  (Mehrtens \cite[p.\;237]{Me93})
\end{quote}
Thus any advocacy of mathematical \emph{modernism} by a marxist
academic must start with defining away \emph{applications} from the
outset as \emph{counter}\-modern--which is precisely what Mehrtens did
in \cite[p.\;523]{Me96} as discussed in Section~\ref{s25}.  For
further details on bellied capitalists see Section~\ref{s47}
(especially note~\ref{f41}).

\section{Mehrtens' book}

In his 1990 book, Mehrtens describes the modernist transformation of
mathematics in the early 20th century in Germany.  According to
Mehrtens, such a transformation is embodied in Hilbert's formalism and
Cantorian set theory.

\subsection{Selective modernist transformation}
\label{s26}

Mehrtens frequently cites Cantor's famous dictum on freedom being the
essence of mathematics, but finesses the issues concerning the reality
of mathematical objects where Cantor held decidedly unmodern views
laced with both theology and metaphysics.

Oddly, Mehrtens does not mention Emmy Noether's school of abstract
algebra, and only briefly notes the impression that van der Waerden's
book \emph{Moderne Algebra} had made on the young Dieudonn\'e.

Whereas Gray treats Klein's \emph{Erlangen Program} (EP) as the first
item of his chapter ``Mathematical Modernism Arrives" in \cite{Gr08},
Mehrtens sees Riemann's \emph{Habilitation} talk of 1854 as the
beginnining of modernism (but see Section~\ref{s311}), and Klein's EP
as the first move of the \emph{Gegenmoderne}.  Having described it in
five pages, he devotes another 19 pages to the interpretation of
footnote III of the EP, titled ``\"Uber den Wert r\"aumlicher
Anschauung,'' i.e., the value of spatial intuition.

Here Mehrtens elaborates the key thesis of his book, namely the divide
between mathematicians of respectively \emph{C-type} and \emph{M-type}
(in our terminology; see Section~\ref{s33}).  Mehrtens quotes, among
others, Heidegger, Foucault, Marx, Kant and Einstein. This list of
authors indicates that the divide has little to do with mathematics
proper, but is concerned rather with the relation between mathematics
and reality.  While the \emph{M-types} (according to Mehrtens) saw
mathematics as a formal system detached from reality (Riemann had
first developed an abstract theory, and then applied it to space), the
\emph{C-types} (Klein among them) insisted that ``there is a true
geometry which is not \ldots{} intended to be merely an illustrative
form of more abstract investigations."

The struggle between \emph{C-types} and \emph{M-types} became more
heated in the decades following, although Mehrtens fails to give a
single example of Klein having actively opposed modernism.
Accordingly, Mehrtens calls Klein an exponent of the \emph{kooperative
Gegenmoderne} \cite[p.\;207]{Me90}.

For Mehrtens, the \emph{Gegenmoderne} culminated in the 1920s and
1930s with Brouwer's attack on Hilbert and with Bieberbach's
\emph{Deutsche Mathematik}, and petered out afterwards.  As is well
known, Bieberbach was hardly mainstream even among the German
mathematicians of his time, as evidenced by the landslide victory of
the counter-proposal to Bieberbach's proposal to adopt a
\emph{F\"uhrerprinzip} within the \emph{Deutsche
Mathematiker-Vereinigung} in the election held on 13\;september 1934.
According to Mehrtens \cite{Me85}, Bieberbach's proposal got 40
no-votes, 11 yes-votes, and 3 empty ballots, whereas the more moderate
counterproposal got 38 yes-votes, 8 no-votes, with 4 abstentions.%
\footnote{The numbers reported by Mehrtens do not quite add up but the
pattern is clear.}

\subsection
{Riemann, Klein, Heidegger, Klein, Riemann}
\label{s311}

Riemann's viewpoint creates a problem for Mehrtens' simplistic
dichotomy of \emph{C-type} versus \emph{M-type}.  As Rowe points out,
\begin{quote}
the contrast [of Hilbert's viewpoint on Euclidean geometry] with
Riemann's viewpoint could not be sharper.  For the latter insisted
that a refined understanding of Euclidean geometry was a dead end.  In
fact, at the very outset of his \emph{Habilitationsvortrag}, Riemann
proclaimed that the study of the foundations of geometry from Euclid
to Legendre--seen as an empirical science--had remained in the dark,
owing to a failure to explore crucial issues or hypotheses concerning
physical measurements.  Natural philosophers thus lacked a general
theory of extended magnitudes or, to use modern language, an
understanding of differential geometry in arbitrary dimensions.
Axiomatics, on the other hand, plays virtually no role in Riemann's
text, least of all speculations about a theory of parallels.
(Rowe \cite[p.\;179]{Ro17})
\end{quote}
Riemann's pervasive influence on modern mathematics ranges from the
theory of Riemann surfaces and Riemannian geometry to the Riemann
hypothesis and the deepest problems in number theory, and reaches as
far as category theory (see Marquis \cite{Ma09}).

Determined to frame Klein as a \emph{C-type}, Mehrtens faces the
formidable challenge of Klein's solid reputation as member of the
Riemannian tradition at G\"ottingen.  Mehrtens' strategy to circumvent
the challenge is a \emph{tour de force} of obfuscation.

Unable to deny Riemann's sterling reputation as a modern, Mehrtens
concedes in \cite[p.\;67]{Me90} that Riemann made `the first move of
modernity'.  However, he then proceeds to paint a different picture,
seeking to portray Riemann as a romantic with pantheistic views.
Mehrtens attempts to back up his picture by collating and
misrepresenting some passages from fragmentary personal notes of
Riemann's (published posthumously by Dedekind and Weber as
\cite{Ri76}) and goes on to claim:
\begin{quote}
In Riemann's philosophical fragments there is talk of an
\emph{earth-soul}, which is a moving, multifaceted
\emph{thought-process}.%
\footnote{In the original: ``In Riemanns philosophischen Fragmenten
ist die Rede von einer ``Erdseele'', die ein bewegter,
vielf{\"a}ltiger ``Denkprozess'' ist.''}
(Mehrtens \cite[p.\;69]{Me90})
\end{quote}
Unable to detect an occurrence of either of the term \emph{Anschauung}
or \emph{Intuition} in these fragments that would provide a figleaf of
respectability for Mehrtens' (partial) C-typing of Riemann, Mehrtens
seeks to connect Riemann to Heidegger, whose name evokes the familiar
sinister associations.

Heidegger's abrupt appearance on page 70 of Mehrtens' book is followed
on page 71 by a quotation of the following difficult passage:
\begin{quote}
Cognition is a kind of representational thinking [\emph{Vorstellen}].
In this \emph{presentation} [\emph{Stellen}] something we encounter
comes to stand [\emph{Stehen}], to a standstill [\emph{Stand}].  What
is encountered and brought to a standstill in representational
thinking is the \emph{object} [\emph{Gegenstand}].%
\footnote{In the original: ``Das Erkennen gilt als eine Art des
Vorstellens.  In \emph{diesem} Stellen kommt etwas, was uns begegnet,
zum Stehen, zum Stand.  Das im Vorstellen zum Stand gebrachte
Begegnende ist der \emph{Gegenstand}'' (Heidegger
\cite[p.\;46]{He57}).  The first sentence does not appear in Mehrtens'
quotation of Heidegger, obscuring the meaning of the already difficult
passage.}
(Heidegger \cite[p.\;23]{He91})
\end{quote}
In his introduction, Mehrtens describes his procedure (including
\emph{Discourse Analysis}) as follows:
\begin{quote}
The sixth chapter, which sums up the interpretatory framework,
and---using \emph{Discourse Analysis}, Semiotics and Semiology, and
Poetology---develops the leitmotifs, systematically and \emph{within}
the three images Eulenspiegel, the Golem, and M{\"u}nchhausen, was the
last chapter that I revised.  (Mehrtens \cite[p.\;9]{Me90}; emphasis
added)
\end{quote}
Examining pages 70--71 in \cite{Me90} and applying such a `Discourse
Analysis' reveals that the names of Riemann (R), Klein (K), and
Heidegger (H) occur here in close succession in the following order:
\begin{quote}
\hfil R, R, K, K, H, K, R. \hfil
\end{quote}
Mehrtens juxtaposes the names and accompanies them by obscure passages
from fragments of Riemann's personal notes and cryptic phrases from
Heidegger taken out of context, skillfully creating a \emph{flou
artistique} that suggests an affinity among the three protagonists,
with an undercurrent of damaging innuendo implicating both Riemann and
Klein in the swamp of Heidegger's well-known forays into politics.%
\footnote{In a similar vein, Mehrtens goes as far as to link both
Poincar{\'e} and Klein with the tainted term \emph{F{\"u}hrer} in
(Mehrtens \cite[p.\;577]{Me90}).  Related innuendo exploiting the
tainted term occurs also on pages 161, 252, 253, 254, 428, 576 in
\cite{Me90}.}

\subsection{Theodor Lessing and race}
\label{s29}

Mehrtens inserts some quotes from an idiosyncratic article by
philosopher Theodor Lessing in the midst of a discussion of Felix
Klein's views in (Mehrtens \cite[pp.\;216--218]{Me90}), including
Lessing's unsourced speculations in \cite[p.\;235]{Le09} regarding
possible Jewish roots of Riemann, Weierstrass, and Klein himself
(Lessing's article, highly critical of Klein's educational reform, was
naturally discussed in Klein's seminar).%
\footnote{Did Dingler's claim that Klein was half-Jewish (see
Section~\ref{s31b}) have Lessing's speculation at its origin?  This
requires further research.}
Lessing's speculations on race in this passage are rather wild-eyed by
modern standards.  Without providing an evaluation of Lessing's
idiosyncratic speculations, Mehrtens abruptly returns to analyzing
Klein's views, thereby implying guilt by association on Klein's part.
A further example of assigning guilt by association is the following
passage:
\begin{quote}
The biologistic-racist discourse, which was an everyday phenomenon in
Wilhelmine Germany, and not only there, and which Klein served and
amplified, was continued in the \emph{Nationalsozialismus}.%
\footnote{\label{f12}In the original: ``Der im wilhelminischen
Deutschland und nicht nur dort allt{\"a}gliche
biologistisch-rassistische Diskurs, den Klein bediente und
verst{\"a}rkte, wurde im Nationalsozialismus fortgef{\"u}hrt.''}
\cite[p.\;217]{Me90}
\end{quote}
The issue of Klein's alleged racism was dealt with in
Section~\ref{s32b}.  For a comparison of Klein's position with that of
Weierstrass see Section~\ref{s210}; for Vahlen see Section~\ref{s310}.
We will deal with the \emph{biologistic} issue separately in
Section~\ref{s44} immediately following.

\subsection
{A nonconventional arsenal: the \emph{biologistischer} stockpile}
\label{s44}

In his book, Mehrtens repeatedly exploits the sinister adjective
\emph{biologistisch}:
\begin{enumerate}
\item{} ``Time and again mathematicians have observed that there is a
generalizable difference between geometric-\emph{anschaulich} and
logical-algebraic thought. \ldots{} Klein conceived of this difference
as a difference of psychic dispositions, and he moreover takes the
step to biologisation when he writes, etc.''%
\footnote{In the original: ``Da{\ss} es jedoch einen
verallgemeinerbaren Unterschied zwischen geometrisch-anschaulich und
logisch-algebraischem Denken gebe, ist von Mathematikern immer wieder
beobachtet worden. \ldots{} Klein fa{\ss}t sie naturalistisch als eine
Verschiedenheit psychischer Dispositionen auf und geht noch den
Schritt zur Biologisierung weiter, wenn er schreibt, etc.''}
\cite[p.\;215]{Me90};
\item{} `biologistic-racist discourse'%
\footnote{In the original: ``biologistisch-rassistische[r] Diskurs.''
This item already appeared in Section~\ref{s29}.}
\cite[p.\;217]{Me90};
\item{} ``The `pure gift of inspiration' links the transcendental
giver with the biologistic conception of `gift' and `selection'.'%
\footnote{In the original: ``Die ``rein geschenkte Eingebung''
verbindet den transzendenten Schenker mit der biologistischen
Konzeption von Begabung und Auslese.''  Note that the German word
\emph{Auslese} (`selection') was routinely used in NS-ideology.}
\\
\cite[p.\;219]{Me90};
\item{} ``The nationalistic or biologistic discourses/narratives of
\emph{Volk} and \emph{Rasse}, as well as the aggressively-heroic
pathos of soldierly masculinity were suitable to be connected with the
self-image of scientists, and also with the [prevailing] idea of
science, this idea having been used by/in countermodernity to defend
meaning and to defend relevance-for-reality.''%
\footnote{In the original: ``Die naturalistischen oder biologistischen
Reden von Volk und Rasse, ebenso wie das aggressiv-heroische Pathos
soldatischer M{\"a}nnlichkeit lie{\ss}sen sich mit dem Selbstbild der
Wissenschaftler und dem Verst{\"a}ndnis von Wissenschaft
verkn{\"u}pfen, mit dem in der Gegegenmoderne Sinn und
Wirklichkeitsbezug der theoretischen Wissenschaft verteidigt worden
waren.''}
\cite[p.\;313]{Me90};
\item{} ``With hindsight it is all too evident that this%
\footnote{Here Mehrtens refers to Klein's comments at Evanston in 1893
(see Section~\ref{s46}); only a single sentence earlier, Klein's name
is mentioned, and associates Klein with ``the imperialist racism of
the dominating `civilized peoples'{}'' (the single quotation marks
around ``civilized peoples' are in the original).}
had the potential to biologistically articulate the possible conflict
among these civilized peoples.''%
\footnote{In the original: ``Da{\ss} darin das Potential lag, den
m{\"o}glichen Konflikt unter jenen 'Kulturnationen' biologistisch zu
artikulieren, ist im nachhinein allzu deutlich.''}
\cite[p.\;345]{Me90};
\item{} ``However, the reform-movement with [its newly coined terms
of] the `modern culture' and the `functional thinking', and with its
biologistic philosophical outlook, does not have anything to do with
the self-definition of science and its fields of application, but is
aimed at gaining legitimacy in the eyes of the general public and [is
also aimed at] the cultural representation of the new realities.''%
\footnote{In the original: ``Die Reformbewegung mit der `modernen
Kultur' und dem `funktionalen Denken', auch mit dem biologistischen
Weltbild, aber hat gerade nicht mit der Selbstdefinition von
Wissenschaft und deren Praxisfeldern zu tun, sondern zielt auf die
{\"o}ffentliche Legitimit{\"a}t und kulturelle Repr{\"a}sentation der
neuen Wirklichkeiten.''}
\cite[p.\;376]{Me90}
\end{enumerate}
Mehrtens' \emph{biologistisches} leitmotif creates a (clearly
intended) impression of continuity between, on the one hand, Mehrtens'
case against Klein and, on the other, Mehrtens' case against the
abuses of the NS era.  Such alleged continuity is spurious (see
further in Section~\ref{s45}).  The kind of role Mehrtens sees for
biology as far as NS ideology is concerned is made crystal clear by
the following passage:
\begin{quote}
A more obvious example from biology is the concept `instinct,' which
was productive in Nazi Germany, both politically and scientifically,
and served to mutually reinforce parts of Nazi ideology and
\emph{biological} ethology\ldots{} (Mehrtens \cite[p.\;229]{Me93};
emphasis added)
\end{quote}
Mehrtens' indiscriminate use of the \emph{biologistischer} stockpile
against Klein is not consistent with standards of meaningful
historical scholarship.

\subsection{Kleinian continuities according to Mehrtens}
\label{s45}

Klein's rather commonplace description of mathematical platonism
appears in his 1926 text:
\begin{quote}
One faction among the mathematicians thinks themselves unrestricted
autocrats in their respective realm, this realm being created of their
own accord, according to their own whim, by logical deductions; the
other faction proceeds from the belief that science pre-exists [the
scientist], being in a state of ideal perfection, and that all that
falls to our lot is to discover, in lucky instants, a limited new
territory, as a piece [of that pre-existing science].%
\footnote{In the original: ``Die eine Gruppe von Mathematikern
h{\"a}lt sich f{\"u}r unbeschr{\"a}nkte Selbstherrscher in ihrem
Gebiet, das sie nach eigener Willk{\"u}r logisch deduzierend aus sich
heraus schaffen; die andere geht von der Auffassung aus, da{\ss} die
Wissenschaft in ideeller Vollendung vorexistiere, und da{\ss} es uns
nur gegeben ist, in gl{\"u}cklichen Augenblicken ein begrenztes
Neuland als St{\"u}ck davon zu entdecken.  Nicht Erfinden nach
Gutd{\"u}nken, sondern Auffinden des ewig Vorhandenen, nicht die
selbstbewu{\ss}te Tat, sondern die vom Bewu{\ss}tsein und Willen
unabh{\"a}ngige, rein geschenkte Eingebung erscheint ihnen als das
Wesen des Schaffens.''}
(Klein \cite[p.\;72]{Kl26})
\end{quote}
Mehrtens reacts to Klein's comments on mathematical Platonism as
follows:
\begin{quote}
It transpires from the evaluative adjectives [employed by Klein] how
Klein places himself; this also agrees with the \emph{`Kleinian
continuities'}, and with the countermodern discourse: the ``logic''
and the ``whim'' of modernity are the \emph{abhorrent} antithesis to
the \emph{rechter} German mathematician, who in pious modesty is
receiving donations of `eternal truth', piecemeal, by way of
``inspiration'' (read: `Anschauung', `Intuition').%
\footnote{In the original: ``Wo er sich selbst einordnet, ist in den
wertenden Beiworten deutlich und f{\"u}gt sich in \emph{Kleins
Kontinuit{\"a}ten} und in den \emph{gegenmodernen Diskurs}: Die
``Logik'' und die ``Willk{\"u}r'' der Moderne sind der abscheuliche
Gegensatz zum \emph{rechten deutschen Mathematiker}, der in frommer
Bescheidenheit die `ewige Wahrheit' st{\"u}ckchenweise aber ``rein''
durch ``Eingebung'' (sprich: Anschauung, Intuition) geschenkt
bekommt'' (emphasis added).}
\cite[p.\;219]{Me90} (emphasis added)
\end{quote}
Mehrtens' conclusion introduces a theological element into the
discussion:
\begin{quote}
The tone [of Klein's voice] is religious. The devil is [incarnated in]
other people; the `fall of man' consists of the [choice to use one's]
whim, and this whim does not bow to the facts of truth and power.%
\footnote{In the original: ``Der Tonfall is religi{\"o}s. Der Teufel,
das sind die anderen; der S{\"u}ndenfall ist die Willk{\"u}r, die vor
den Gegebenheiten der Wahrheit und der Macht den Nacken nicht
beugt.''}
(ibid.)
\end{quote}
The discrepancy between Klein's comment and Mehrtens' intemperate
reaction to it requires no amplification.  The Talmudic dictum ``kol
haposel, bemumo posel'' (whoever disqualifies others, [it is] in his
own blemish [that he] disqualifies [them]) provides an insight into
this passage, where \emph{Mehrtens} himself coins the expression
\emph{wertende Beiworte}: Klein never used the adjective
\emph{abscheulich} [`abhorrent'] with regard to any of the ideas
mentioned: `logic', `whim/arbitrariness', or `modernity'.%
\footnote{Nor did Klein refer to non-formalist mathematicians as the
`right' mathematicians.  Mehrtens wishes to see a political statement
in places where there is none.  In keeping with an apparently
irrepressible urge to paint Klein as a right-wing precursor of worse
things to come, Mehrtens plays on an ambiguity of the German adjective
\emph{recht} (it could mean either right-minded, righteous, or
right-wing).  Furthermore, Klein did not use the religious term
\emph{fromm} anywhere, nor \emph{ewig}, nor \emph{Teufel}, nor
\emph{S{\"u}ndenfall}, nor \emph{Macht}.  Mehrtens' claim to the
contrary is sheer fabrication.  It is not Klein but Mehrtens who
struck a religious tone, contrasting with the merely old-fashioned
tone of Klein's comment.}

The \emph{continuity} alleged in Mehrtens' passage is all Mehrtens',
not Klein's.  Mehrtens' procedure here is another example of
\emph{massaging the evidence}.

\subsection{The butterfly model from Klein to Vahlen}
\label{s310}

Writes Mehrtens:
\begin{quote}
Theodor Vahlen, who after 1933 was an executive official in the
ministry, and a professor in Berlin, gave, on [15 may] 1923, an
address on assuming his office as rector [of University of
Greifswald], in the usual tradition, speaking on `Value and Essence'
of mathematics, and cited \emph{Klein's racist distinctions} within
mathematics, and sharpened them into open antisemitism.%
\footnote{In the orginal: ``Theodor Vahlen, nach 1933 leitender
Ministerialbeamter und Professor in Berlin, hielt 1923 eine
Rektoratsrede in einschl\"agiger Tradition `Wert und Wesen' der
Mathematik und zitierte \emph{Kleins rassistische Unterscheidungen} in
der Mathematik und versch{\"a}rfte sie zu offenem Antisemitismus''
(emphasis added).  Mehrtens' reference to ``Professor in Berlin, hielt
1923 eine Rektoratsrede, etc.'' may give the impression that Vahlen
gave the antisemitic address at Berlin, since the word `Greifswald' is
relegated to a footnote, and `Universit{\"a}tsreden' is omitted,
making it appear more significant than it actually was (Greifswald is
a small city).}
\cite[p.\;310--311]{Me90} (emphasis added)
\end{quote}
What we object to most is Mehrtens' judgmental phrase `Klein's racist
distinctions.'  Here Mehrtens again equivocates on the meaning of a
loaded term (see Section~\ref{s32b}), and furthermore implies a
continuity (and perhaps even organic necessity) between Klein's
remarks on ethnic differences on the one hand, and Vahlen's
``sharpened'' open antisemitism, on the other.  In the background of
Mehrtens' remarks is an unspoken endorsement of the \emph{butterfly
model} of the evolution of ideas, contrasted with the \emph{Latin
model} (see Ian Hacking \cite[p.\;119]{Ha14}).

Mehrtens' allegation of continuity between the views of Klein
and those of Vahlen (who eventually reached the rank of
\emph{SS-Brigadef\"uhrer})%
\footnote{See St\"ower \cite[p.\;145, note\;349]{St12}.}
is an additional instance of assigning guilt by association, as well
as a character smear against Klein.

In line with the Talmudic dictum cited in Section~\ref{s45}, Mehrtens
seeks to argue a case of \emph{biological} inevitability between the
cocoon of Klein's remarks on intuition and race, on the one hand, and
the dragon-butterfly (see e.g., Benisch \cite{Be91}) of the
\emph{Nationalsozialismus}, on the other.

\subsection{Imperialist fight and a 1908 cartoon}
\label{s47}

Mehrtens implicates Klein in no less than an ``imperialist competitive
fight:''
\begin{quote}
In the same year [1901, in a speech] before the `support club'%
\footnote{This is a reference to an association charged with carrying
out the reforms in the teaching of mathematics and science in
German-speaking high school teaching.}
Klein emphasizes ``the importance of [higher] education in mathematics
and the sciences in the competitive fight of the nations, and of the
general cultural significance of high school education, and of a
capable class of high school teachers supported by public trust.''  If
the higher totality is culture, then the concrete totality is the
nation in the \emph{imperialist competitive fight.}
%
%
\ldots{} The goal was to define a \emph{modern} culture in such a way,
and to make it a state affair in such a way, that the professionals in
natural- and social sciences, economics and bureaucracy rise to the
rank of a state elite.%
\footnote{In the original: ``Im gleichen Jahr betonte Klein vor dem
F{\"o}rderverein $\rangle\rangle$die Bedeutung des
mathematisch-naturwissenschaftlichen Studiums im Konkurrenzkampf der
Nationen und die allgemeine Kulturbedeutung der h{\"o}heren Schule und
eines leistungsf{\"a}higen, vom {\"o}ffentlichen Vertrauen getragenen
Standes der Studienr{\"a}te.$\langle\langle$
Wenn das h{\"o}here Ganze die Kultur ist, dann ist das konkrete Ganze
die Nation im imperialistischen Konkurrenzkampf. 
%
%
\ldots{} Es kam darauf an, eine \emph{moderne} Kultur so zu definieren
und zur Staatsangelegenheit zu machen, da{\ss} die beteiligten
Professionellen aus Natur- und Sozialwissenschaft, Wirtschaft und
B{\"u}rokratie gemeinsam in den Rang einer staatlichen Elite
aufr{\"u}ckten.''}
\cite[p.\;361]{Me90} (emphasis on `modern' in
the original; emphasis on `imperialist competitive fight' added)
\end{quote}
The tone of the above passage is consistent with the identification of
the source of the enmity toward modernism on the part of marxist
academics and Hilton Kramer's diagnosis thereof (see
Section~\ref{s26b}).

That the driving force behind Mehrtens' animus toward Klein is marxist
detestation of capitalism is confirmed by Mehrtens' analysis of a
cartoon reproduced at \cite[p.\;381]{Me90}.  The cartoon represents
the activities of the \emph{G{\"o}ttinger Vereinigung} founded by
Klein, formed to encourage interactions between the Academy and
Industry.  The cartoon depicts recognizably professorial types
exchanging books for proverbial moneybags brought by recognizably
capitalist types.%
\footnote{\label{f41}Mehrtens describes the cartoon as follows: ``A
contemporary caricature is more precise [than the logo] in that regard
(Figure 5).  Therein, at the signpost pointing towards the
`G{\"o}ttinger Vereinigung', professors, clad in cap and gown, and
with a large and a small book tucked under their arms, meet with the
capitalists, unmistakable with their bellies and top hats, who carry a
large and a small moneybag.  After exchanging the small book for the
small moneybag, they push along in pairs, engrossed in conversation''
(translation ours).  Mehrtens' sarcastic tone indicates that he views
the cartoon as a satirical caricature; see note~\ref{f23}.}
Hovering over the scene are a sun-like Klein and an angelic-looking
Althoff (see Section~\ref{s22b}), depicted as bestowing blessings upon
the congregants.  The cartoon is dated at 1908 in also \cite{He80},
\cite{HS91}, \cite{Oh87}.

Mehrtens' purpose in reproducing the cartoon seems to be to compare it
to the logo of the \emph{G{\"o}ttinger Vereinigung} which depicts a
allegorical scene alluding to academic-industrial interactions.  Thus,
Mehrtens declares the cartoon to be a `more accurate' (``genauer'')
representation of the activities of the \emph{G{\"o}ttinger
Vereinigung} than the official logo (see \cite[p.\;382]{Me90}).%
\footnote{\label{f23}While Mehrtens interprets the cartoon as a
satirical caricature of the goings-on at the \emph{G{\"o}ttinger
Vereinigung}, Hermann and Sch\"onbeck note in \cite[Caption on
p.\;355]{HS91} that the cartoon was sent out (as part of the
invitation letter) by the organizers of the conference themselves
(apparently in self-deprecating humor).  Ohse et
al.\;\cite[p.\;373]{Oh87} note that the cartoon displays the following
text: ``Gruss vom Festkommers zur Feier des 10-j{\"a}hrigen Bestehens
der G{\"o}ttinger Vereinigung G{\"o}ttingen, 22. Febr. 1908''
(translation: ``Greetings from the \emph{Commercium} on the occasion
of the 10th anniversary of the existence of the G{\"o}ttinger
Vereinigung G{\"o}ttingen, 22 february 1908'').  The source for the
cartoon is given as follows: Archiv der Aerodynamischen
Versuchsanstalt, Prandtl Zimmer, Ordner $\rangle\rangle$G{\"o}ttinger
Vereinigung 1905-1919$\langle\langle$.}
Mehrtens' own attitude toward such activities is spelled out clearly
enough in the pages leading up to the cartoon:
\begin{quote}
The disposition over technical knowledge is being regulated via legal
norms, by patent laws, and by laws for the protection of various
secrets, in such a way that said knowledge is \emph{co-opted} for use
within the capitalist market, and for use within the \emph{military
state.}%
\footnote{In the original: ``Mit Patentrecht und Geheimnisschutz wird
die Verf{\"u}gung {\"u}ber technisches Wissen durch Rechtsnormen so
reguliert, da{\ss} es f{\"u}r den kapitalistischen Markt und den
\emph{Milit{\"a}rstaat zugerichtet} wird'' (emphasis added).}
\cite[pp.\;378--379]{Me90} (emphasis added)
\end{quote}
Mehrtens' skillful insinuation of a connection between Klein and the
\emph{military state} finds its full expression in Mehrtens' treatment
of Theodor Vahlen; see Section~\ref{s310}.

\subsection{Modernism in America}

According to Mehrtens, after the war modernism triumphed and
culminated in the works of Nicolas Bourbaki \cite[p.\;320]{Me90}.
However, this is a simplistic view of mathematical modernism.  Thus,
Mehrtens makes no mention of the debate over modernist ideas in
America.  Notably, Marshall Stone claimed that
\begin{quote}
while several important changes have taken place in our conception of
mathematics or in our points of view concerning it, the one which
truly involves a revolution in ideas is the discovery that mathematics
is entirely independent of the physical world (except that thinking
takes place in the brain).  (Stone \cite[p.\;716]{St61})
\end{quote}
Stone's idea of mathematics as independent of the physical world was
eventually picked up by Quinn; see Bair et al.\;\cite{18b}.
%
%

In response to Stone's paper a symposium was held.  Here Richard
Courant, who was in many respects Klein's heir, had this to say:
\begin{quote}
Certainly mathematical thought operates by abstraction; mathematical
ideas are in need of abstract progressive refinement, axiomatization,
crystallization. It is true indeed that important simplification
becomes possible when a higher plateau of structural insight is
reached.  \ldots{} Yet, the life blood of our science rises through
its \emph{roots}; these roots reach down in endless ramification deep
into what might be called reality, whether this ``reality" is
mechanics, physics, biological form, economic behavior, geodesy, or,
for that matter, other mathematical substance already in the realm of
the familiar.  (Courant in Carrier et al.\;\cite{Ca62}; emphasis
added)
\end{quote}
Courant's \emph{tree} metaphor, later elaborated by Kline, contrasted
with the Bourbaki--Tucker \emph{city} metaphor; see Phillips
\cite{Ph14} for an analysis.  Courant continued:
\begin{quote}
Abstraction and generalization is not more vital for mathematics than
individuality of phenomena and, before all, not more than inductive
intuition. Only the interplay between these forces and their synthesis
can keep mathematics alive and prevent its drying out into a dead
skeleton. \ldots{} We must not accept the old blasphemous nonsense
that the ultimate justification of mathematical science is `the glory
of the human mind'.  (Courant in\;\cite{Ca62})
\end{quote}
Similarly, a letter by Morris Kline and others claimed that the
relevance of mathematics for students emerges from concrete
situations, not formalisms; Kline was to argue his point most
forcefully in \cite{Kl73}.  The history of modernism in mathematics is
more complex than Herbert Mehrtens is willing to grant.

\subsection{Klein's Evanston lectures}
\label{s46}

Klein gave twelve lectures at Northwestern University at Evanston in
1893.  The published version of the lectures occupies 98 pages.  His
sixth lecture contains the following passage:

\begin{quote}
[Alfred] K\"opcke,%
\footnote{Mathematician, 1852--1927, doctoral student of Leo
K\"onigsberger 1875 at Heidelberg.}
of Hamburg, has advanced the idea that our space-intuition is exact as
far as it goes, but so limited as to make it impossible for us to
picture to ourselves curves without tangents.

On one point Pasch does not agree with me, and that is as to the exact
value of the axioms.  He believes -- and this is the traditional view
-- that it is possible finally to discard intuition entirely, basing
the whole science on the axioms alone. I am of the opinion that,
certainly, for the purposes of research it is always necessary to
combine the intuition with the axioms.  I do not believe, for
instance, that it would have been possible to derive the results
discussed in my former lectures, the splendid researches of Lie, the
continuity of the shape of algebraic curves and surfaces, or the most
general forms of triangles, without the constant use of geometrical
intuition.

Pasch's idea of building up the science purely on the basis of the
axioms has since been carried still farther by Peano, in his logical
calculus.

Finally, it must be said that the degree of exactness of the intuition
of space may be different in different individuals, perhaps even in
different races. It would seem as if a strong naive space-intuition
were an attribute pre-eminently of the \emph{Teutonic race}, while the
critical, purely logical sense is more fully developed in the
\emph{Latin} and \emph{Hebrew races}.  A full investigation of this
subject, somewhat on the lines suggested by Francis Galton in his
researches on heredity, might be interesting.  (Klein
\cite[p.\;45--46]{Kl94}; emphasis added)
\end{quote}
Klein's tentative comments on Teutonic, Hebrew, and Latin races,
occupying a total of 6 lines in a 98-page text, were made in the
context of
\begin{itemize}
\item
a discussion of Kopcke's views concerning space-intuition;
\item
an appreciation of Pasch's accomplishment vis-a-vis the foundations;
\item
Klein's disagreement with Pasch with regard to the role of intuition
in the creative process.
\end{itemize}
We propose our own tentative interpretation of the ``Teutonic, Latin,
Hebrew'' (T, L, H) comment.  Confronted with the striking diffence in
mathematical style between K\"opcke (T) of Hamburg and Pasch~(H),
Klein is led to speculate on possible ethnic origins of such
differences, with the remark about the Latin races added because of
Klein's awareness of the axiomatic/foundational contribution of the
Italian mathematician Peano\;(L).  Historian Sanford Segal reacted to
Klein's statements as follows:
\begin{quote}
For Klein, though certainly a conservative nationalist, was also
certainly no anti-Semite--he had helped bring first Hermann Minkowski
and then Edmund Landau as well as Karl Schwarzschild to G\"ottingen,
and spoke favorably of the emancipation of the Jews in Prussia in 1812
(by Napoleon Bonaparte): ``With this action a large new reservoir of
mathematical talent was opened up for our country, the powers of
which, coupled with the increase therein achieved by the French
emigrants, very soon proved itself fruitful in our science.''  He also
had good relations with a number of Jewish mathematicians and had an
extensive correspondence with his good friend Max Noether (Emmy's
father): eighty-nine letters from Noether and 129 from Klein, often
addressed ``Lieber Noether''.  (Segal \cite[p.\;270]{Se03})
\end{quote}
We will analyze Mehrtens' rather different reaction in
Section~\ref{s49} immediately following.

\subsection{Mehrtens on Klein's 6-line comment}
\label{s49}

Historian Sanford Segal's relaction to Klein's statements appears
above (Section~\ref{s46}).  The reaction to the same statements by
Mehrtens is on record:
\begin{quote}
It was less antisemitism or racist nationalism which found its
expression in [Klein's] statements, rather it was the imperialist
racism of the dominating civilized peoples; this means that Klein gave
an internationalist argumentation, by implicitly mentioning the limits
of the Internationalism of his time.  That [Klein's comments] had the
potential to enable others to articulate possible conflict between
those `civilized nations' biologistically, is, with hindsight, all too
clear.%
\footnote{In the original: ``Weniger Antisemitismus oder rassistischer
Nationalismus kam hier zum Ausdruck als eher der imperialistische
Rassismus der dominierenden `Kulturv{\"o}lker'; das hei{\ss}t, Klein
argumentierte internationalistisch, indem er implizit die Grenzen des
zeitgen{\"o}ssischen Internationalismus mit ansprach.  Da{\ss} darin
das Potential lag, den m{\"o}glichen Konflikt unter jenen
`Kulturnationen' biologistisch zu artikulieren, ist im nachhinein
allzu deutlich'' (an excerpt of this quotation already appeared in
Section~\ref{s44}).}
(Mehrtens \cite[p.\;345]{Me90})
\end{quote}
It appears that \emph{antisemitism, racism, imperialism, domination,}
and (especially the sinister) \emph{biologism} (see Section~\ref{s44})
are the first things that come to Mehrtens' mind when reading Klein's
tentative comments on space-intuition and possible differences among
ethnicities in intellectual outlook.  As we already noted in
Section~\ref{s32b}, Mehrtens' \emph{racism} slur against Klein is
based on equivocation.

Mehrtens' diatribe against Klein rings particularly hollow to a reader
aware of Mehrtens' cynical exploitation of the infamous \emph{yellow
star} badge in his comments on Hilbert:
\begin{quote}
Hilbert's [set-theoretic] paradise is a dictatorship\ldots{} The
yellow star is, when viewed mathematically, pure set-building.%
\footnote{In the original: ``Das Hilbertsche Paradies ist eine
Diktatur\ldots{} Der gelbe Stern ist, mathematisch gesehen, reine
Mengenbildung.''}
(Mehrtens \cite[p.\;460]{Me90})
\end{quote}
Mehrtens' flippant ``set-building'' comment comes disturbingly close
to an odious trivialisation of the \emph{yellow star} badge, a tragic
symbol of the Holocaust.

\subsection{Mehrtens and Gray on Hilbert}
\label{s315}

Mehrtens presents Hilbert as the arch-modernist (``Generaldirektor der
Moderne'').  In Gray's words:
\begin{quote}
Hilbert in particular has a major role, and his work is presented
ironclad as a program to make all of mathematics abstract, axiomatic
and internally self-consistent.  (Gray \cite[p.\;10]{Gr08})
\end{quote}
However, such a view of Hilbert's posture is based on a very selective
reading of his works, taking into account only the \emph{Grundlagen
der Geometrie} and his work on the logical foundations of mathematics
in the 1920s.

On Hilbert's \emph{Zahlbericht}, Corry writes: ``The general idea of
algebra as the discipline dealing with algebraic structures is still
absent from Hilbert's work on algebraic number theory'' (Corry
\cite[p.\;154]{Co04}).

The same can be said concerning Hilbert's papers on integral
equations, which were written between 1904 and 1910 in the language of
classical analysis, and are far from the foundational work on abstract
function spaces that was done at the same time by Fr\'echet, Riesz and
Schmidt.  Hilbert never ventured into a study of such spaces.

Corry quotes a passage from a 1905 lecture, which illustrates
Hilbert's rather \emph{in}formal approach to axiomatics:
\begin{quote}
The edifice of science is not raised like a dwelling, in which the
foundations are first firmly laid and only then one proceeds to
construct and to enlarge the rooms.  Science prefers to secure as soon
as possible comfortable spaces to wander around and only subsequently,
when signs appear here and there that the loose foundations are not
able to sustain the expansion of the rooms, it sets about supporting
and fortifying them.  This is not a weakness, but rather the right and
healthy path of development.  (Hilbert quoted in Corry \cite{Co07})
\end{quote}
The following quotation from the foreword to the translation of
Hilbert and Cohn-Vossen's \emph{Anschauliche [!] Geometrie} \cite{HC}
shows Hilbert in agreement with Klein with regard to the interplay
between abstraction and intuition:
%
%
\begin{quote}
In mathematics, as in any scientific research, we find two tendencies
present. On the one hand, the tendency toward abstraction seeks to
crystallize the logical relations inherent in the maze of material
that is being studied, and to correlate the material in a systematic
and orderly manner. On the other hand, the tendency toward intuitive
understanding fosters a more immediate grasp of the objects one
studies, a live rapport with them, so to speak, which stresses the
concrete meaning of their relations.

As to geometry, in particular, the abstract tendency has here led to
the magnificent systematic theories of Algebraic Geometry, of
Riemannian Geometry, and of Topology; these theories make extensive
use of abstract reasoning and symbolic calculation in the sense of
algebra.  Notwithstanding this, it is still as true today as it ever
was that intuitive understanding plays a major role in geometry. And
such concrete intuition is of great value not only for the research
worker, but also for anyone who wishes to study and appreciate the
results of research in geometry.
\end{quote}
Here Hilbert can hardly be said to be involved in an ironclad program
of making all of mathematics abstract and axiomatic.  The empiricist
inclinations in Hilbert's work in geometry are studied by Corry
\cite{Co06}.  Notes Rowe:
\begin{quote}
As Michael Toepell has convincingly shown, \emph{Anschauung} held an
important place in Hilbert's geometrical work, and he was by no means
convinced that one could ultimately dispense with it altogether.
(Rowe \cite[p.\;197]{Ro94})
\end{quote}
As far as the relationship between Klein and Hilbert is concerned,
Mehrtens claims that Klein, described as the jupiter of G\"ottingen,
``tolerated Hilbert's formalism'':
\begin{quote}
It pertains to the [concept of] \emph{Unordnung} [disorder] that Felix
Klein, the Jupiter of G{\"o}ttingen, and the defender of
\emph{Anschauung}, became an advocate of university education for
women and tolerated Hilbert's formalism, and is not so disorderly
after all, if one takes into consideration Liebermann, the Secession,
and the male self-construction in the portrait [of Klein by
M. Liebermann] (cf. 7.2).%
\footnote{In the original: ``Da{\ss} Felix Klein, der Jupiter
G{\"o}ttingens und Verteidiger der Anschauung, f{\"u}r das
Frauenstudium eintrat und den Formalismus Hilberts tolerierte,
geh{\"o}rt zur Unordnung und ist im {\"u}brigen so ordnungslos nicht,
denkt man an Liebermann, die Secession und die m{\"a}nnliche
Selbstkonstruktion im Portr{\"a}t (vgl.\;7.2).''}
(Mehrtens \cite[p.\;577]{Me90})
\end{quote}
Since the bulk of Hilbert's work with Bernays on Formalism and
foundations started at about the time of Klein's death, Mehrtens'
comment is an \emph{ahistorical collage} in addition to involving an
unhelpful stereotype.

\section{Conclusion: Mehrtens' tools}

In sum, this particular marxist historian has exploited a variety of
tools in his analysis of Klein that ranges from massaging the evidence
(Sections~\ref{s16}, \ref{s45}) and character smear
(Sections~\ref{s32b}, \ref{s310}) to assigning guilt by association
(Sections~\ref{s29}, \ref{s310}) and ahistorical collage
(Section~\ref{s315}).

On occasion, Mehrtens sheds the mask of a marxist historian and
engages in what is discernibly bourgeois yellow journalism, as in
Mehrtens' tale concerning the dress embroidered with images of
analytic curves with which Felix Klein allegedly ``covered the body''
of his bride, the tale in question being further embroidered by
Mehrtens' tasteless comments at \cite[p.\;214]{Me90}.

Mehrtens' claim that Hilbert was interested in production rather than
meaning applies only to the claim's author.  Mehrtens' portrayal of
Klein as countermodern is contrary to much historical evidence and
must be rejected.

\section*{Acknowledgments}

We are grateful to Simcha Horowitz for helpful suggestions.  The
influence of Hilton Kramer (1928--2012) is obvious.

\end{document}